\documentclass[10pt]{article}

\usepackage[utf8]{inputenc}
\usepackage{amsmath} 
\usepackage{amssymb} 
\usepackage{amsthm}
\usepackage{xcolor}
\usepackage{parskip}
\usepackage[T1]{fontenc}

\newtheorem{proposition}{Proposition}[subsection]
\newtheorem{Remarque}[proposition]{Remark}
\newtheorem{conjecture}[proposition]{Conjecture}
\newtheorem{fait}[proposition]{Fact}
\newtheorem{definition}[proposition]{Definition}
\newtheorem{theoreme}[proposition]{Theorem}

\newtheorem{corollaire}[proposition]{Corollary}
\def\Ind#1#2{#1\setbox0=\hbox{$#1x$}\kern\wd0\hbox to
0pt{\hss$#1\mid$\hss}
\lower.9\ht0\hbox to 0pt{\hss$#1\smile$\hss}\kern\wd0}

\def\Notind#1#2{#1\setbox0=\hbox{$#1x$}\kern\wd0\hbox to 0pt{\mathchardef
\nn="3236\hss$#1\nn$\kern1.4\wd0\hss}\hbox to 0pt{\hss$#1\mid$\hss}\lower.9\ht0
\hbox to 0pt{\hss$#1\smile$\hss}\kern\wd0}

\title{Ranked definably linear quasi-Frobenius groups
\footnote{ keywords : Frobenius groups, groups of finite Morley tank, involutions
MSC 2020 : 20N05, 20F11, 03C60}}

\author{Samuel Zamour \footnote{Univ Paris Est Creteil, Univ Gustave Eiffel, CNRS, LAMA UMR8050, F-94010 Creteil, France, e-mail : samuel.zamour@u-pec.fr}}
\date{11/06/2024}
\begin{document}
\maketitle
\begin{abstract}
We describe the classification of ranked definably linear quasi-Frobenius groups of odd type : dihedral configurations are isomorphic to $PGL(2,K)$ for $K$ an algebraically closed field of characteristic other than two; Frobenius groups are spilt and solvable if the characteristic of the underlying field is positive. 
\end{abstract}

\section{Introduction}
The classical groups $\mbox{PGL}(2,\mathbb{C})$, $\mbox{SO}(3,\mathbb{R})$ and $\mbox{GA}(1,\mathbb{C})$ have a common group-theoretic structure: there is a proper subgroup which has finite index in its normalizer and trivial intersection with its distinct conjugates. This common structure is closely related to the incidence geometry induced by involutions. 
\\
These group-theoretical configurations were first considered in \cite{DW}; they generalize the well-known notion of \textit{Frobenius groups} - the proper subgroup is self-normalizing - which appears notably in the classification of finite simple groups. In our paper \cite{Zam1}, we proposed to call these configurations \textit{quasi-Frobenius groups}, and we pursued their study in the context of model theory. More specifically, we would like to classify ranked quasi-Frobenius groups with involutions. Roughly speaking, ranked groups are abstract groups with a dimension function on definable sets generalizing the Zariski dimension for algebraic varieties; finite groups and algebraic groups over algebraically closed fields are typical examples (see below for a precise definition). Some new results on ranked quasi-Frobenius groups were obtained by Alt\ inel, Deloro and Corredor in \cite{ACD}; in particular, the so-called "odd degree" configurations are eliminated. Notice also that ranked Frobenius groups are the objet of a renowned interest and have been recently studied in \cite{CT} and \cite{PoiF}.
\\
A \textit{ranked group} (or group of finite Morley rank) is an abstract group with a dimension function $RM$ on definable/interpretable sets satisfying the following axioms: 
\begin{itemize}
    \item Axiom A (monotonicity): Let $A$ be a definable set. Then $RM(A)\geq n+1$ if and only if $A$ contains an infinite number of pairwise disjoint definable subsets of rank greater than or equal to $n$.
    \item Axiom B (definability) : If $f: A\rightarrow B$ is a definable function between two definable sets $A$ and $B$, then for any integer $n$, the subset $\{ b\in B : RM(f^{-1}(b))=n \}$ is definable.
    \item Axiom C (additivity): If $f:A\rightarrow B$ is a surjective definable function, and if $RM(f^{-1}(b))=n$ for every $b\in B$, then $RM(A)=n+RM(B)$.
    \item Axiom D (elimination of the infinite quantifier): For any definable function $f:A\rightarrow B$, there exists an integer $m$ such that for any $b\in B$, the fiber $f^{-1}(b)$ is infinite as soon as it contains at least $m$ elements.
\end{itemize}
The Cherlin-Zilber conjecture asserts that a ranked infinite simple group is isomorphic to an algebraic group. For a general reference on ranked groups, we refer to \cite{BN}.
\\
Before stating the main conjecture about quasi-Frobenius groups, let us give a precise definition.
\begin{definition}
    Let $C<G$ be such that $[N_G(C):C]<\infty$ and $C\cap C^g=\{1\}$ for all $g\notin N_G(C)$. If $G$ and $C$ are ranked connected groups, we say that $C<G$ is a ranked quasi-Frobenius group. We say it is of odd type if it contains involutions and if it is $U_2^{\perp}$. In addition, there are two cases depending on the value of $[N_G(C):C]$ \cite[Theorem B]{ACD} : 
    \begin{itemize}
        \item If $[N_G(C):C]=1$, we speak of a ranked Frobenius group.
        \item (dihedral configuration) If $[N_G(C):C]$ is even, we speak of a ranked quasi-Frobenius group of even degree.
        
    \end{itemize}
\end{definition}

\begin{conjecture} \cite[A1 conjecture]{DW}
Let $C<G$ be a ranked quasi-Frobenius group of odd type. 
    \begin{itemize}
        \item If $[N_G(C):C]=2$, then $G\simeq PGL(2,K)$ for $K$ an algebraically closed field of characteristic other than 2.
        \item If $[N_G(C):C]=1$, then $G$ is split and solvable. 
        \end{itemize}
\end{conjecture}
We will show that the conjecture is verified for a subclass of ranked groups: ranked definably linear groups. We say that $G$ is a \textit{ranked definably linear group} (we will also say a \textit{definably linear group over a ranked field}) if $G$ is definable in a (potentially enriched) ranked field and if it is isomorphic to a subgroup of a general linear group over this field.
\\
We begin with the dihedral case. Then, we move to Frobenius groups (section 2). Finally, in a third section, we prove the following theorem :
\begin{theoreme} \label{standard Frobenius quasigroups} Let $G=A\rtimes C$ be a split ranked quasi-Frobenius group of odd type. Then $C\leq \mbox{GL}_n(K)$, for a field $K$. Moreover, one of the following cases occurs :
\begin{enumerate}
    \item $K$ is a field of positive characteristic $p>0$, and $C$ is a good torus.
    \item $K$ is a field of characteristic zero and if $C$ is not abelian, then $C$ has a structure analogous to that of a bad group: it is covered by its maximal linear tori.
\end{enumerate}
\end{theoreme} We will deduce as a corollary that if $C<G$ is a ranked definably linear Frobenius group of odd type over a field of positive characteristic, then $G$ is solvable.
\subsection*{Acknowledgements}
This paper is based on the author's Phd thesis directed by Frank Wagner in Université Lyon 1. We would like to thank him for his guidance and for his numerous and helpful comments and remarks. We are also very grateful to Adrien Deloro for his corrections and advises.
\section{The dihedral case}
First, we recall some definitions and results concerning ranked quasi-Frobenius groups. 

\begin{fait}\label{quasi-groupe de Frobenius}
    \begin{enumerate}
        \item \cite[Proposition 1]{DW} If $C<G$ is a ranked quasi-Frobenius group of even degree and of odd type. Then $[N_G(C):C]=2$ and $C$ is abelian inverted by an involution.
        \item \cite[Proposition 4.2.3]{Zam1} Let $C<G$ be a ranked quasi-Frobenius of odd type. Suppose that there exists a definable connected normal subgroup $A\lhd G$ such that $A\cap C=\{1\}$. Then there exists a definable connected abelian subgroup $A_1$ such that $G=A_1\rtimes C$
        \item \cite[Theorem 4.5.4]{Zam1} Let $C<G$ be a  ranked quasi-Frobenius group of even degree and odd type. There exists a weakly standard Borel subgroup, i.e., a non-nilpotent Borel subgroup which intersects $C$ in a infinite subset.
        \item \cite[Theorem 1.0.6]{Zam1} Let $C<G$ be a  ranked quasi-Frobenius group of even degree and odd type. Assume that there exists a \textit{strongly standard Borel}, i.e there exists a non-nilpotent Borel subgroup strictly containing $C$. Then $G\simeq PGL(2, K)$, for $K$ an algebraically closed field of characteristic not two.
    \end{enumerate}
\end{fait}

\begin{theoreme}
Let $C<G$ be a  definably linear quasi-Frobenius group of even degree and of odd type over a ranked field. Then $G\simeq \mbox{PGL}_2(K)$, for $K$ an algebraically closed field of characteristic not two.
\end{theoreme}
\begin{proof}
Recall that $C$ is abelian (Fact \ref{quasi-groupe de Frobenius}). Assume that $C$ is a Borel subgroup. Let $B$ be a weakly standard (non-nilpotent) Borel subgroup. By conjugation of the Borel subgroups \cite[Proposition 2.11]{Mus}, there exists $g\in G$ such that $B=C^g$, hence $B$ is nilpotent, a contradiction. 
\\
Consequently, there exists a strongly standard Borel subgroup. We conclude with Fact \ref{quasi-groupe de Frobenius}. \end{proof}

\section{The Frobenius group case}
We can now pass to the classification of ranked definably linear Frobenius groups of odd type.
\\
The situation of algebraic Frobenius groups is fully elucidated by the following fact:
\begin{fait}\cite[Theorem 4]{Her}\label{algebraic Frobenius}
Let $C < G$ be an algebraic Frobenius group. Then $G$ is split.
\end{fait}

Here is a useful fact about 2-torsion :
\begin{fait}\label{2-central twist} \cite[Theorem B]{DW}
Let $G$ be a ranked connected group, $U_2^{\perp}$, such that the Borel subgroups are nilpotent and hereditarily conjugate. Then the (potentially trivial) 2-torsion is central.
\end{fait}
\begin{theoreme}
 Let $C<G$ be a ranked connected definably linear Frobenius group of odd type. Then $G$ is split.\end{theoreme}
\begin{proof} 
We may assume that $R(G)=\{ 1 \}$: otherwise, $G$ would contain a normal abelian subgroup $A$ such that $A\cap C=\{ 1 \}$ and by \cite[Lemma 11.21]{BN}, $G$ would be split.
If the characteristic of the underlying field is positive, then by \cite[Théorème 2.6]{Mus}, $G$ is isomorphic to a finite product $S_1\times...\times S_r $ of simple algebraic groups $S_i$. As all involutions are contained in $\bigcup_G C^g$, $S_i\cap C\neq \{1\}$, for $1\leq i \leq r$. Therefore, we may suppose that $G=S_1$  and $G$ is a simple algebraic group, a contradiction.
\\
If the characterictic is zero, by \cite[Théorème 2.7]{Mus}, $G$ is a finite product $S_1\times...\times S_r\times C$, where each $S_i$ is a simple algebraic group and $C$ is a definable subgroup without unipotents. Therefore, if $G$ did not contain only semisimple elements, we could assume that $(C\cap S_1)<S_1\neq \{1\}$ is a simple algebraic Frobenius group, a contradiction. As a consequence, $G$ contains only semisimple elements. Let $B$ be a Borel subgroup of $G$; as $B'$ is unipotent, we have $B'=\{1\}$ and so $B$ is abelian. By \cite[Proposition 2.11]{Mus}, the Borel subgroups are hereditarily conjugate, and by Fact \ref{2-central twist}, we conclude that $Z(G)\neq \{1\}$, a contradiction. 
\end{proof}
\section{The split case}
The properties of ranked definably linear groups allow us to describe the structure of split ranked connected quasi-Frobenius groups of odd type. In a certain sense, these groups are "standard": we can plunge the complement into a definably linear group over a ranked field and in positive characteristic $p>0$, it is a good torus. We obtain a result analogous to the main theorem of \cite{CGNV}, which describes the structure of linear groups (over an algebraically closed field) that are sharply transitive for their natural action.
\\
Before proving the main theorem of this section, let us recall an easy consequence of \cite[Théorème 2.6]{Mus}.
\begin{fait}\label{borel nilpotents}
Let $G$ be a ranked connected definably linear group over a field of positive characteristic. If the Borel subgroups of $G$ are nilpotent, then $G$ is nilpotent.    
\end{fait}
\begin{proof}
If $G$ is solvable, then it is its own Borel subgroup, hence nilpotent.
\\
Suppose that $G$ is not solvable. If all proper definable connected subgroups are nilpotent, then $G$ is a bad group. Otherwise, there is a proper definable connected non-solvable subgroup whose Borel subgroups are nilpotent (the Borel subgroups of $G$ are nilpotent). By the descending chain condition on definable sugbroups, we may assume that $G/R(G)$ is a bad group. However, by \cite[Théorème 2.6]{Mus}, $G/R(G)$ is isomorphic to a finite product of simple algebraic groups whose Borel subgroups are therefore nilpotent, a contradiction \cite[Proposition B. § 21.4]{Hum}. 
\end{proof}
\begin{theoreme} \label{standard Frobenius quasigroups} Let $G=A\rtimes C$ be a split ranked quasi-Frobenius group of odd type. Then $C\leq \mbox{GL}_n(K)$, for a field $K$. Moreover, one of the following cases occurs :
\begin{enumerate}
    \item $K$ is a field of positive characteristic $p>0$, and $C$ is a good torus.
    \item $K$ is a field of characteristic zero and if $C$ is not abelian, then $C$ has a structure analogous to that of a bad group: it is covered by its maximal linear tori.
\end{enumerate}
\end{theoreme}
\begin{proof}
First of all, notice that the involution contained in $C$ inverts $A$ which is therefore abelian ($C_A(i)$ is indeed finite).
\\
Let $A_1\leq A$ be an infinite definable $C$-minimal subgroup. The group $A_1$ is abelian and connected. The action of $C$ remains free hence faithful : let $1\neq c\in C$ be such that $A_1\leq C_G(c)$, we have $A_1\leq C_G(c)^{\circ}\leq C$, a contradiction. In particular, $C$ is an infinite group of automorphisms of $A_1$. There are two cases: either $A_1$ is divisible without torsion or $A_1$ has a finite exponent.
\\
a) $A_1$ is of finite exponent $p\neq 0$: Let $B\leq C$ be a maximal definable connected solvable subgroup, and $A_2\leq A_1$ be a (infinite) $B$-minimal definable subgroup. Note that the action remains faithful. By \cite[Theorem 9.8]{BN}, the subgroup $B$ is abelian and by \cite[Theorem 9.1]{BN}, we can interpret a field $K$ of characteristic $p>0$ such that $B\leq K^{\times}$. By \cite{Wag4}, $B$ is a good torus; therefore the Borel subgroups of $C$ are abelian and (hereditarily) conjugate. By Fact \ref{2-central twist}, the infinite 2-torsion of $C$ is central, in particular, $Z(C)$ is infinite. By \cite[Theorem 9.5]{BN}, there exists a definable field $K$ of characteristic $p>0$, such that $A_1$ is a finite-dimensional $K$-vector space and $C$ acts linearly. By Fact \ref{borel nilpotents}, we conclude that $C$ is nilpotent, hence abelian.
\\
b) $A_1$ is divisible without torsion. By \cite{LW}, there exists a definable field $K$ (of characteristic zero) such that $C\leq \mbox{GL}_n (K)$. The subgroup $C$ is therefore definably linear over a ranked field. Moreover, $C$ contains only semisimple elements : let $1\neq c$ be an unipotent element; as the only eigenvalue is one, the associated eigenspace $\{a\in A_1 : c.a=a\}$ is of dimension at least one, and therefore, $c$ has a infinite number of fixed points, a contradiction. As a consequence, the Borel subgroups of $C$ are abelian. We can also assume that $C$ is not solvable. 
\end{proof}
\begin{corollaire}
Let $C < G$ be a definably linear Frobenius group of odd type over a ranked field of characteristic $p>0$. Then $G$ is solvable.
\end{corollaire}
\begin{Remarque}
    Independently, Poizat establishes some results on Frobenius groups in \cite{PoiF}. In particular, he shows that pseudo-locally finite Frobenius groups are split.
\end{Remarque}
\bibliography{definably linear}
\bibliographystyle{plain}
\end{document}